\documentclass{amsart}
\usepackage{amssymb,amsmath,amscd}
\usepackage{amsthm}
\usepackage{amsbsy}
\usepackage{enumerate}
\usepackage[all]{xy}

\def\ord{\text{ord}}
    \def\nz{\hbox{\text{0}{ \raise 3pt \hbox{\kern -10pt
    \vrule width 8pt height 0.4 pt }\kern 2pt}}}

\begin{document}
\title{REES VALUATIONS}
\author{By Shreeram S. Abhyankar and William J. Heinzer}
\address{Mathematics Department, Purdue University, West Lafayette, IN 47907, USA.}
\email{ram@cs.purdue.edu, heinzer@math.purdue.edu}
\begin{abstract}
We study dicritical divisors and Rees valuations. 

2000 Mathematics Subject Classification: Primary 14A05.

Key Words: contact number, dicritical divisor, tangent cone.
\end{abstract}
\date{}
\maketitle

\markboth{}{}

\baselineskip 20 pt

{{\bf Section 1: Introduction.}} 
In Example (8.3) of our recent paper \cite{AH2}, we considered the local ring
$R$ of the point $(0,0,0)$ on the irreducible surface
$$
Z^m-F_1(X,Y)\dots F_h(X,Y)=0
$$ 
where $m>h>0$ are integers such that $m$ is nondivisible by the 
characteristic of the ground field $K$, and $F_1(X,Y),\dots,F_h(X,Y)$ are 
pairwise coprime homogeneous linear polynomials. 
Now geometrically speaking, a section of the tangent cone of the surface at 
$(0,0,0)$ consists of the $h$ lines $f_1(X,Y)=0,\dots,f_n(X,Y)=0$
in the $(X,Y)$-plane passing through $(0,0)$, and they ``give rise to'' all the
distinct $h$ elements of $\mathfrak D(R,M)$, where $M$ is the maximal ideal
$M(R)$ of $R$ and  $\mathfrak D(R,M)$ is the set of all dicritical divisors
of $M$ in $R$.
See Section 2 for the definition of dicritical divisors and their relation
to Rees valuations. The analytical (topological) theory of dicritical divisors 
was developed in \cite{Art,EiN,Fou,LeW,MaM}. 
It was algebracized in \cite{Ab9} to \cite{Ab14} and \cite{AH1,AH2,AbL}.  

In (3.2) and (3.3) of Section 3 we shall show that $M$ is a normal ideal in $R$, 
i.e., for every positive integer $p$, the ideal $M^p$ is a complete ideal in 
the normal local domain $R$ and, if $F_0(X,Y)$ is a homogeneous linear 
polynomial which is coprime to $F_1(X,Y),\dots,F_h(X,Y)$ then the ideal 
generated by $F_0(X,Y)^p$ and $Z^p$ is a reduction of $M^p$; as a special 
example we can take $F_1(X,Y)\dots F_h(X,Y)=X^h-Y^h$ and $F_0(X,Y)^p=Y^p$.
Actually, in (3.2) and (3.3), we shall prove a generalized higher dimensional 
version of this result by two different methods.
In (3.4) to (3.7) we shall make some relevant observations.
In (4.9) to (4.11) we shall raise some questions.

In the above phrase ``give rise to,''  we are trying to mimic Max Noether's 
concept of ``infinitely near in the first neighborhood.'' Thus, 
making a QDT = Quadratic Transformation centered at the point $(0,0,0)$ of the
above displayed surface, i.e., substituting $(X,Y)=(ZX',ZY')$ and
factoring out the ``exceptional locus'' $Z^h$ we get
the transformed surface
$$
Z^{m-h}-F_1(X',Y')\dots F_h(X',Y')=0.
$$ 
Let $V_j$ be the local ring of the line $Z=F_j(X',Y')=0$ on the transformed
surface. Then $V_1,\dots,V_h$ are exactly all the distinct members of
$\mathfrak D(R,M)$.

In Section 4 we shall elucidate the contact number $c(R,V,W)$ which appears
in Proposition (3.5) of \cite{Ab12}. Here $V$ and $W$ are prime divisors of a
two dimensional regular local domain $R$.
We shall prove some commutativity properties of the contact number and relate
them to calculations of the local intersection multiplicity of two curves
at a simple point of an algebraic or arithmetical surface. 

\centerline{}

{{\bf Section 2: Terminology and Preliminaries.}} 
We shall use the notation and terminology introduced in Sections 2 to 9 of 
\cite{AH1} and Sections 2 to 9 of \cite{AH2} which themselves were based on 
\cite{Ab8} to \cite{Ab14}. Relevant background material can be found in 
\cite{Ab1} to \cite{Ab7} and \cite{Nag,NoR,Zar}.

\centerline{}

Referring to pages 145-161 of \cite{Ab8}, 
for the foundations of models, recall that: 

$\mathfrak V(A)=\{A_P:P\in\text{spec}(A)\}=$ the modelic spec of a domain $A$.

$S^{\mathfrak N}=$ the set of all members of $\mathfrak V(\overline S)$ which
dominate $S$ where $\overline S$ is the integral closure of a quasilocal domain
$S$ in its quotient field QF$(S)$.

$U^{\mathfrak N}=\cup_{B\in U}B^{\mathfrak N}$ for any set $U$ of quasilocal
domains.

$\mathfrak W(A,J)=\cup_{0\ne x\in J}\mathfrak V(A[Jx^{-1}])=$ the modelic 
blowup of $A$ at a nonzero ideal $J$ in a domain $A$; note that
$Jx^{-1}=\{y/x:y\in J\}$; see pages 152-160 of \cite{Ab8}.

$\overline D(L/A)=$ the set of all valuation rings $V$ with quotient field $L$
such that $A\subset V$ where $A$ is a subring of $L$. For this and other related
definitions see Section 2 of \cite{Ab11}. In the cited reference we had defined
the completion of an ideal $J$ in a domain $A$ with quotient 
field $L$ only when $A$ is normal. i.e., integrally closed in $L$; now we extend 
this without assuming $A$ to be normal by saying that the {\bf completion}
$\overline J$ of $J$ is always obtained by putting
$\overline J=\cap_{V\in\overline D(L/A)}(JV)$; note that $\overline J$ is clearly
an ideal in the integral closure of $A$ in $L$.

As usual $\mathbb N$ (resp: $\mathbb N_+$) denotes the set of all nonnegative 
(resp positive) integers. The set of all nonzero elements in a ring $A$ is denoted
by $A^\times$. 

 For any set $U$ of quasilocal domains and any $i\in\mathbb N$,
$U_i$ denotes the set of all $i$-dimensional members of  $U$.

$\mathfrak W(S,J)^\Delta_i=$ the set of all $i$-dimensional members of
$\mathfrak W(S,J)$ which dominates $S$ where $J$ is a nonzero ideal in a quasilocal
domain $S$ and $i\in\mathbb N$.

$\mathfrak D(S,J)=(\mathfrak W(S,J)_1^\Delta)^{\mathfrak N}=$ the 
{\bf dicritical set} of a nonzero ideal $J$ in a quasilocal domain $S$; members of
this set are called {\bf dicritical divisors} of $J$ in $S$. 

\centerline{}

More generally:

$\mathfrak W(S,J,M)^\Delta_i=$ the set of all $i$-dimensional members $T$ of
$\mathfrak W(S,J)$ such that $M\subset M(T)$ where
$J\subset M$ are nonzero ideals in a domain $S$ and $i\in\mathbb N$.

$\mathfrak W^*(S,J,M)^\Delta_i=$ the set of all $i$-dimensional members $T$ of
$\mathfrak W(S,J)^{\mathfrak N}$ such that $M\subset M(T)$ where
$J\subset M$ are nonzero ideals in a domain $S$ and $i\in\mathbb N$.

$\mathfrak D(S,J,M)=(\mathfrak W(S,J,M)^{\Delta}_1)^{\mathfrak N}=$ 
the {\bf dicritical set} of a pair of nonzero ideals $(J,M)$ in a domain $S$
with $J\subset M$; members of this set are again
called the {\bf dicritical divisors} of $(J,M)$ in $S$. 

$\mathfrak D^*(S,J,M)=(\mathfrak W(S,J,M)^{\mathfrak N})^\Delta_1=$ the 
{\bf starred dicritical set} of a pair of nonzero ideals $(J,M)$ in a domain $S$
with $J\subset M$; members of this set are 
called the {\bf starred dicritical divisors} of $(J,M)$ in $S$. 

\centerline{}

We make the following observations concerning dicritical divisors.

(I) If $J\subset M$ are nonzero ideals in a domain $S$ and 
$N$ is an ideal in $S$ with $J\subset N$ such that rad$_SN=\text{rad}_SM$
then for all $i\in\mathbb N$ we have 
$\mathfrak W(S,J,N)^\Delta_i=\mathfrak W(S,J,M)^\Delta_i$ with
$\mathfrak W^*(S,J,N)^\Delta_i=\mathfrak W^*(S,J,M)^\Delta_i$, and we have
$\mathfrak D(S,J,N)=\mathfrak D(S,J,M)$
with $\mathfrak D^*(S,J,N)=\mathfrak D^*(S,J,M)$. If $M=S$ then for all
$i\in\mathbb N$ we have $\mathfrak W(S,J,M)^\Delta_i=\emptyset$ with
$\mathfrak W^*(S,J,M)^\Delta_i=\emptyset$, and we have
$\mathfrak D(S,J,M)=\emptyset$ with
$\mathfrak D^*(S,J,M)=\emptyset$.

(II) If $J\subset M$ are nonzero ideals in a quasilocal domain $S$ with $M=M(S)$ 
then $\mathfrak W(S,J,M)^\Delta_i=\mathfrak W(S,J)^\Delta_i$ for all 
$i\in\mathbb N$, and we have $\mathfrak D(S,J,M)=\mathfrak D(S,J)$.

(III) If $R$ is a noetherian domain and $I\subset M$ are nonzero nonunit ideals 
in $R$ then $\mathfrak D(R,I,M)$ is a 
finite set of DVRs. As in $(5.6)(\dagger^*)$ of \cite{Ab8},
this follows from (33.2) on page 115 of \cite{Nag} or (33.10) on page 118 of 
\cite{Nag}; here (33.2) is called the Krull-Akizuki Theorem and, according to 
page 218 of \cite{Nag}, (33.10) may be called the Mori-Nagata Theorem.
Note that if $M\subset\text{rad}_RI$ then $\mathfrak D(R,I,M)\ne\emptyset$.

(IV) Let us note that, for a nonzero ideal $J$ in a noetherian domain $R$, 
members of $\mathfrak D^*(R,J,J)$ are called 
{\bf Rees valuation rings} of $J$ in $R$, and their order functions are called
{\bf Rees valuations} of $J$ in $R$. Note that, by the Mori-Nagata Theorem 
(33.10) on page 118 of \cite{Nag}, Rees valuation rings are DVRs.

\centerline{}

{{\bf Section 3: Normality of the Maximal Ideal.}} 
In (3.1) and (3.2) we use the set-up of a QDT; 
for details see pages 534-577 of \cite{Ab8}
which are based on pages 7-45, 155-192, 262-283 of \cite{Ab3}.
Let us start by discussing:  

\centerline{}

(3.1) TANGENT CONES AND QUADRATIC TRANSFORMATIONS. Let $d$ be a positive
integer and let $X_1,\dots,X_{d+1}$ be indeterminates over a field $K$.
Let $m\ge h$ be positive integers. Let
$$
G=G(X_1,\dots,X_{d+1})=\sum_{h\le i\le m}G_i\in C=K[X_1,\dots,X_{d+1}]
$$
be irreducible with
$$
G_i=G_i(X_1,\dots,X_{d+1})\in C
$$
such that
depending on $i$, $G_i=0$ or $G_i$ is a homogeneous polynomial of degree $i$,
and we have $G_h\ne 0\ne G_m$.
Note that then $G_h=0$ is the tangent cone of the hypersurface
$G=0$ at the origin $(0,\dots,0)$. Let $R$ be the local ring of the
origin on the hypersurface and let $M$ be the maximal ideal $M(R)$ in $R$.
Let $B$ be the affine coordinate ring of the hypersurface. Then
$B=C/(GC)=K[x_1,\dots,x_{d+1}]$ where
$x_1,\dots,x_{d+1}$ are the images of $X_1,\dots,X_{d+1}$ in $C/(GC)$ under
the residue class epimorphism $\Phi:C\to C/(GC)$ and we
have identified $K$ with its image in $C/(GC)$. The $d$-dimentional local domain
$R$ is the localization of $B$ at the maximal ideal $(x_1,\dots,x_{d+1})B$
and we have $M=(x_1,\dots,x_{d+1})R$. The proof of the following RELATION 
between the tangent cone and the modelic blowup $\mathfrak W(R,M)$ is essentially
contained in pages 534-577 of \cite{Ab8}.

\centerline{}

(I) RELATION.
{\it Assume that the hyperplane $X_{d+1}=0$ is not a component of the tangent cone,
i.e., $G_h(X_1,\dots,X_d,0)\ne 0$, and let 
$A_R=R[x_1/x_{d+1},\dots,x_d/x_{d+1}]$. Then we have
$\mathfrak W(R,M)^\Delta\subset\mathfrak V(A_R)$. Hence in particular
$\mathfrak D(R,M)=(\mathfrak V(A_R)_1)^{\mathfrak N}$.}

\centerline{}

Now let us prove the following REDUCTION property of the tangent cone.

\centerline{}

(II) REDUCTION. {\it Let $H_1,\dots,H_d$ be homogeneous linear members of
$C$ which are linearly independent over $K$ and which do not divide $G_h$.
Then the ideal in $R$ generated by $\Phi(H_1),\dots,\Phi(H_d)$
is a reduction of $M$.}

\centerline{}

PROOF. Without loss of generality we may assume that
$H_1=X_1,\dots,H_d=X_d$. Now the assumption that
$H_1,\dots,H_d$ do not divide $G_h$ is equivalent to saying that
$G_h(0,\dots,0,X_{d+1})\ne 0$. The form ring $F_{(R,M)}(M)$ is, in a natural
manner, isomorphic to the associated graded ring grad$(R,M)$ as discussed on
pages 272-277 of \cite{Ab8} and hence, in view of the material on these pages,
by (6.1) of \cite{AH1} we see that the ideal $(x_1,\dots,x_d)R$ is a reduction
of $M$.

\centerline{}

(3.2) FIRST (HIGH-SCHOOL) METHOD. Let $d,m,h$ be a positive integers with
$m>h$, let $X_1,\dots,X_{d},Z$ be indeterminates over a field $K$, and let 
$$
G=G(X_1,\dots,X_{d},Z)=Z^m+\sum_{h\le i\le m-1}G_i\in C=K[X_1,\dots,X_{d},Z]
$$
be irreducible with
$$
G_i=G_i(X_1,\dots,X_{d},Z)\in C
$$
such that
depending on $i$, $G_i=0$ or $G_i$ is a homogeneous polynomial of degree $i$.
Let $R$ be the local ring of the origin $(0,\dots,0)$ on the hypersurface
$G=0$, and let $M$ be the maximal ideal $M(R)$ in $R$.
Assume $(\sharp)$ and $(\sharp\sharp)$ sandwiched between items (17) and (18) of the
following DD = Detailed Description. In DD we shall show that
$\mathfrak D(R,M)$ consists of $h$ distinct DVRs $V_1,\dots,V_h$ and upon
letting (32) of DD we have (33), (45), and (46) of DD.

\centerline{}

DETAILED DESCRIPTION.  Consider the polynomial rings
$$
C'=K[X_1,\dots,X_d]\subset K[X_1,\dots,X_d,Z]=C
\leqno(1)
$$
where $d$ is a positive integer and 
$X_1,\dots,X_d,Z$ are indeterminates over a field $K$. Recalling that the
mspec of a ring is the set of all maximal ideals in it, let
$$
Q_{C'}=(X_1,\dots,X_d)C'\in\text{mspec}(C')
$$
with 
$$
R_{C'}=C'_{Q_{C'}} \;\;\text{ and }\;\;L_{C'}=\text{QF}(C')=K(X_1,\dots,X_d)
$$
and
$$
Q_C=(X_1,\dots,X_d,Z)C\in\text{mspec}(C)
$$
with 
$$
R_C=C_{Q_C} \;\;\text{ and }\;\;L_C=\text{QF}(C)=K(X_1,\dots,X_d,Z).
$$
Note that then $R_C$ is a $(d+1)$-dimensional regular local domain which 
dominates the $d$-dimensional regular local domain $R_{C'}$.  Let
$$
x_1=X_1/Z,\dots,x_d=X_d/Z
\leqno(2)
$$
and note that then $x_1,\dots,x_d,Z$ may be regarded as indeterminates over $K$.
Let
$$
A_{C'}=K[x_1,\dots,x_d]\subset K[x_1,\dots,x_d,Z]=A_C.
\leqno(3)
$$
Then
$$
\begin{cases}
$$
\text{$C\subset A_C$ are $(d+1)$-variable polynomial rings}\\
\text{with a common quotient field $L_C$}\\
\text{while $C'\not\subset A_{C'}\not\subset C'$ are 
$d$-variable polynomial rings}\\
\text{which are subrings of the field $L_C$.}\\
\end{cases}
$$

Let 
$$
\text{$m=t+h$ where $t$ and $h$ are positive integers.}
\leqno(4)
$$
Let 
$$
G=G(X_1,\dots,X_d,Z)=\sum_{h\le i\le m}G_i\in C\setminus C'
\leqno(5)
$$
with
$$
G_i=G_i(X_1,\dots,X_d,Z)\in C
\leqno(6)
$$
such that
$$
\text{depending on $i$, $G_i=0$ or $G_i$ is homogeneous of degree $i$}
\leqno(7)
$$
and we have
$$
G_h\ne 0\ne G_m.
\leqno(8)
$$
Since $G\in C\setminus C'$, we get $C'\cap(GC)=0$. Therefore we can construct 
an overring $B$ of $C'$ together with 
$$
\begin{cases}
\text{$C'$-epimorphism $\Phi:C\to B=K[X_1,\dots,X_d,z]$}\\
\text{ with $\Phi(Z)=z$ and $\text{ker}(\Phi)=GC$}\\
\end{cases}
\leqno(9)
$$
which may be depicted by the following commutative diagram

\centerline{}

\begin{equation*}
\CD
GC@>{\text{ker}(\Phi)=GC}>{\text{inj}}>C=K[X_1,\ldots,X_d,Z]
@>{\Phi}>{\text{sur}}>K[X_1,\ldots,X_d,z]=B\\
@.@A{\text{inj}}AA     @AA{=}A \\
@.C'=K[X_1,\ldots,X_d]@>{\text{inj}}>> K[X_1,\ldots,X_d,z]=B
\endCD
\end{equation*}

\centerline{}

where inj and sur indicate injective and surjective maps respectively.

Henceforth
$$
\text{assume that $G$ is irreducible}
\leqno(\dagger)
$$
i.e., equivalently, assume that $B$ is a domain. Let
$$
Q=(X_1,\dots,X_d,z)B\in\text{mspec}(B) 
\leqno(10)
$$ 
with 
$$
R=B_{Q} \;\;\text{ and }\;\;L=\text{QF}(B)=K(X_1,\dots,X_d,z)
\leqno(11)
$$
Note that then $R$ is a $d$-dimensional local domain which dominates $R_{C'}$ and 
$\Phi$ extends to a unique $(R_{C'})$-epimorphism
$$
\Phi^*:R_C\to R\;\;\text{ with }\;\;\text{ker}(\Phi^*)=GR_C.
$$
Henceforth
$$
\text{assume that }G_m(X_1,\dots,X_d,Z)=Z^m
\leqno(\ddagger)
$$
and let
$$
g=g(x_1,\dots,x_d,Z)=\frac{G(Zx_1,\dots,Zx_d,Z)}{Z^h}
$$
and
$$
x'_1=X_1/z,\dots,x'_d=X_d/z
\leqno(12)
$$
and
$$
A'=K[x'_1,\dots,x'_d]
\;\;\text{ with }\;\;L'=K(x'_1,\dots,x'_d)
\leqno(13)
$$
and
$$
A=K[x'_1,\dots,x'_d,z]
\;\;\text{ with }\;\;L_A=K(x'_1,\dots,x'_d,z).
\leqno(14)
$$
Then clearly
$$
g=g(x_1,\dots,x_d,Z)=Z^t+\sum_{1\le i\le t}g_iZ^{t-i}\in A_C
\leqno(15)
$$
where
$$
g_i=g_i(x_1,\dots,x_d)=\frac{G_{m-i}(Zx_1,\dots,Zx_d,Z)}{Z^{m-i}}\in A_{C'}
\leqno(16)
$$
with
$$
g_t\ne 0
$$
and
$$
\begin{cases}
\text{$A'$ is a $d$-variable polynomial ring with QF$(A')=L'$,}\\
\text{$A$ is a $d$-dimentional noetherian domain with QF$(A)=L=L_A$,}\\
\text{and we have the subring inclusions $A'\subset B$ and $B\subset A$}\\
\end{cases}
$$
and $\Phi$ can be uniquely extended to 
$$
\begin{cases}
\text{$K$-epimorphism $\phi:A_C\to A$}\\
\text{ with $\phi(x_1,\dots,x_d,Z)=x'_1,\dots,x'_d,z$
 and $\text{ker}(\phi)=gA_C$}\\
\end{cases}
\leqno(17)
$$
which may be depicted by the following commutative diagram

\centerline{}

\begin{equation*}
\CD
@.C=K[X_1,\ldots,X_d,Z]@.K[x'_1,\ldots,x'_d]=A'\\
@.@V{\text{inj}}VV     @VV{\text{inj}}V \\
gA_C@>{\text{ker}(\phi)=gA_C}>{\text{inj}}>A_C=K[x_1,\ldots,x_d,Z]
@>{\phi}>{\text{sur}}>K[x'_1,\ldots,x'_d,z]=A\\
@.@A{\text{inj}}AA     @AA{=}A \\
@.A_{C'}=K[x_1,\ldots,x_d]@>{\text{inj}}>> K[x'_1,\ldots,x'_d,z]=A\\
@.@.     @A{X_i=zx'_i}A{\text{inj}}A \\
@.@.K[X_1,\ldots,X_d,z]=B\\
\endCD
\end{equation*}

\centerline{}

where inj and sur indicate injective and surjective maps respectively.

 For a moment suppress assumption $(\dagger)$ but henceforth assume that
$$
\sum_{h\le i<m}G_i\in FC\text{ with }F=F(X_1,\dots,X_d)=\prod_{1\le j\le h}F_j\\
\leqno(\sharp)
$$
where
$$
\begin{cases}
F_j=F_j(X_1,\dots,X_d)\in C'=K[X_1,\dots,X_d]\\
\text{and $F_1,\dots,F_h$ are nonzero homogeneous linear polynomials}\\
\text{which are coprime;}\\
\text{here coprime means the ideals $F_1C',\dots,F_hC'$ are distinct.}
\end{cases}
\leqno(\sharp\sharp)
$$
By Eisenstein's Criterion we see that 
$$
(\sharp)+(\sharp\sharp)\Rightarrow(\dagger).
\leqno(18)
$$ 

Upon letting
$$
f=f(x_1,\dots,x_d)=\prod_{1\le j\le h}f_j
$$
with 
$$
f_j=f_j(x_1,\dots,x_d)=\frac{F_j(Zx_1,\dots,Zx_d)}{Z}
$$
we see that
$$
f_j=f_j(x_1,\dots,x_d)\in A_{C'}
$$
are such that
$$
\begin{cases}
\text{$f_1,\dots,f_h$ are coprime nonzero homogeneous linear polynomials}\\
\text{where coprime means the ideals $f_1A_{C'},\dots,f_hA_{C'}$ are distinct}\\
\end{cases}
$$
and we have
$$
g_{t}A_{C'}=fA_{C'}
\;\;\text{ and }\;\;
\sum_{1\le i\le t}g_iZ^{t-i}\in fA_C.
$$

Now upon letting
$$
f'=f'(x'_1,\dots,x'_d)=f(x'_1,\dots,x'_d)\in A'
\leqno(19)
$$
and 
$$
f'_j=f'_j(x'_1,\dots,x'_d)=f_j(x'_1,\dots,x'_d)\in A'
\;\;\text{  for }\;\;1\le j\le h
\leqno(20)
$$
we see that
$$
f'=f'(x'_1,\dots,x'_d)=\prod_{1\le j\le h}f'_j
\leqno(21)
$$
and
$$
\begin{cases}
\text{$f'_1,\dots,f'_h$ are coprime nonzero homogeneous linear polynomials}\\
\text{where coprime means the ideals $f'_1A',\dots,f'_hA'$ are distinct}\\
\end{cases}
\leqno(22)
$$
and 
$$
\phi(g_t)A'=f'A'\;\;\text{ and }\;\;\phi(g_i)\in f'A'
\;\;\text{  for }\;\;1\le i\le t.
\leqno(23)
$$

 Upon letting
$$
\text{$V'_j=$ the localization $A'_{f'_jA'}$
\;\;for \;\;$1\le j\le h$} 
\leqno(24)
$$
we see that
$$
\text{$V'_1,\dots,V'_h$ are distinct DVRs with quotient field $L'$.} 
\leqno(25)
$$
Now $L=L'(z)$ and $z$ is a root of
$$
Z^t+\sum_{1\le i\le t}\phi(g_i)Z^{t-i}
$$
and hence upon letting
$$
\text{$V_j=$ the integral closure of $V'_j$ in $L$}
\;\;\text{ for }\;\;1\le j\le h 
\leqno(26)
$$
by (12) to (25) we see that
$$
\text{$V_1,\dots,V_h$ are distinct DVRs with quotient field $L$} 
\leqno(27)
$$
and for $1\le j\le h$ we have
$$ 
V_j(y)=\begin{cases}
t+1&\text{if }y=F_j\\
1&\text{if }y=z\\
1&\text{if $y=F_i$ with }i\in\{1,\dots,h\}\setminus\{j\}\\
e&\text{if $y\in C'\setminus(F_jC')$ is homogeneous of degree $e\in\mathbb N$}\\
\end{cases}
\leqno(28)
$$
where the third line is a special case of the fourth line.
Upon letting
$$
M=M(R)
\leqno(29)
$$
by (10) and (11) we see that
$$
M=(X_1,\dots,X_d,z)R.
\leqno(30)
$$
In view of (3.1)(I), by (12) to (30) we see that
$$
\text{$V_1,\dots,V_h$ are exactly all the distinct members of 
$\mathfrak D(R,M)$.}
\leqno(31)
$$

Upon letting
$$
I_p=(M(V_1)^p\cap R)\cap\dots\cap(M(V_h)^p)\cap R)
\;\;\text{ for all }\;\;p\in\mathbb N 
\leqno(32)
$$
WE CLAIM THAT
$$
I_p=M^p\;\;\text{ for all }\;\;p\in\mathbb N. 
\leqno(33)
$$
Namely, by (28) to (32) we get $M^p\subset I_p$. Given any $y\in I_p$ we shall
show that $y\in M^p$ and this will prove (33). Suppose if possible that
$y\not\in M^p$. 

Then, for some $0\le q<p$, we must have 
$y\in M^q$ with $y\not\in M^{q+1}$. Now clearly
$$ 
y=y'+\sum_{0\le i\le q}a_i(X_1,\dots,X_d)z^{q-1}
\;\;\text{ with }\;\;y'\in M^{q+1}
\leqno(34)
$$
where 
$$
\text{$a_i(X_1,\dots,X_d)\in C'$ is either $0$ or homogeneous of degree $i$}
\leqno(35)
$$
Let
$$
\Lambda=\{i:0\le i\le q\text{ and }a_i(X_1,\dots,X_d)\not\in FC'\}.
\leqno(36)
$$
By (4) to (9), $(\ddagger)$, $(\sharp)$, $(\sharp\sharp)$ we see that
$0\ne F\in C'$ is homogeneous of degree $h$ with $F\in M^{h+1}$;
consequently by (28) and (34) we see that if $i\in\{0,\dots,q\}\setminus\Lambda$
then $a_i(X_1,\dots,X_d)z^{q-i}\in M^{q+1}$; therefore upon letting
$$
y''=\sum_{i\in\Lambda}a_i(X_1,\dots,X_d)z^{q-i}
\leqno(37)
$$
by (34) we see that
$$
y-y''\in M^{q+1}.
\leqno(38)
$$
Since $y\not\in M^{q+1}$, by (37) and (38) we get
$$
\Lambda\ne\emptyset.
\leqno(39)
$$
Upon letting
$$
y^*=y''/z^q
\leqno(40)
$$
by (12), (35), and (37) we see that
$$
y^*=\sum_{i\in\Lambda}a_i(x'_1,\dots,x'_d)
\leqno(41)
$$
where, for each $i\in\Lambda$, 
$$
\text{$a_i(x'_1,\dots,x'_d)\in A'\setminus(f'A')$ is homogeneous of degree $i$}
\leqno(42)
$$
and hence by (21) and (22) we see that for some $j\in\{1,\dots,h\}$ we have
$$
a_i(x'_1,\dots,x'_d)\in A'\setminus(f'_jA')
\;\;\text{ for all }\;\;i\in\Lambda 
\leqno(43)
$$
and now by (41) to (43) we conclude that
$$
y^*\in A'\setminus(f'_jA')
$$
and therefore by (26), (28), (38), (39), and (40) we get
$$
V_j(y)=q
$$
wnich contradicts the assumption $y\in I_p$. Therefore we must have 
$y\in M^p$. This completes the proof of (33).

\centerline{}

(44) OBSERVE THAT item (28) can be sharpened by saying that, for $1\le j\le h$,
the DVR $V_j\cap L_{C'}$ can be described thus:

We can write $C'=$ the polynomial ring $K[F_j,Y_{2j},\dots,Y_{dj}]$ 
where $Y_{2j},\dots,Y_{dj}$ are suitable homogeneous linear members of $C'$.
Taking a $(t+1)$-th root $Y_{1j}$ of $F_j$ we may regard $C'$ as a subring of 
the polynomial ring $C_j=K[Y_{1j},\dots,Y_{dj}]$.
We get a DVR $\overline V_j$ with quotient field $K(Y_{1j},\dots,Y_{dj})$
such that for every $0\ne y\in C_j$ we have $\overline V_j(y)=\ord(y)$ 
where $\ord(y)=$ the degree of the smallest degree term when we express $y$ as 
a polynomial in $Y_{1j},\dots,Y_{dj}$. Now clearly 
$V_j\cap L_{C'}=\overline V_j\cap L_{C'}$.

\centerline{}

Taking $H_d=Z$ in (3.1)(II) and invoking (3.4)(II) below we see that
$$
\begin{cases}
\text{if there exist homogeneous linear members $H_1,\dots,H_{d-1}$ of $C'$}\\
\text{which are linearly independent over $K$}\\
\text{and which do not divide $F_1\dots F_h$ in $C'$ then, 
for every $p\in\mathbb N$,}\\
\text{the ideal $(H_1^p,\dots,H_{d-1}^p,z^p)R$ is a reduction of $M^p$.}
\end{cases}
\leqno(45)
$$

\centerline{}

Finally, by standard arguments we see that
$$
\begin{cases}
\text{if $G_i=0$ for $h<i<m$}\\
\text{and $m$ is nondivisibe by the characteristic of $K$}\\
\text{then $R$ is normal}.
\end{cases}
\leqno(46)
$$

\centerline{}

(3.3) SECOND (REES-RING) METHOD.
Let $I$ be an ideal in a nonnull ring $R$. Referring to pages
272-277 of \cite{Ab8} for definitions, for the associated graded ring 
$$
\text{grad}(R,I)=\bigoplus_{n\in\mathbb N}(I^n/I^{n+1})
$$
we have the following.

\centerline{}

(3.3.1) LEMMA. Assume that the ring grad$(R,I)$ is reduced, i.e., has no nozero 
nilpotent element. Then for every positive integer $c$,
the ideal $I^c$ coincides with its integral closure in $R$.

\centerline{}

PROOF. Let if possible $c$ be a positive integer such that $I^c$ does not coincide
with its integral closure in $R$. Then there exists $x\in R\setminus I^c$ such that
$x$ is integral over $I^c$. Let
$$
x^n+\sum_{1\le i\le n}a_ix^{n-i}=0\;\;\text{ with }\;\; a_i\in I^{ci}
$$ 
be an equation of integral dependence. Clearly there is a unique nonegative integer
$r<c$ such that $x\in I^r\setminus I^{r+1}$. Let
$$
\overline x=\text{lefo}_{(R,I)}(x)=\text{the leading form of $x$ relative to }(R,I).
$$
Then $0\ne\overline x\in\text{grad}(R,I)$. Now for all $i$ we have
$ci+(n-i)r\ge nr+1$ and hence $a_ix^{n-i}\in I^{nr+1}$. Therefore
$$
x^n=-\sum_{1\le i\le n}a_ix^{n-i}\in I^{nr+1}
$$ 
and hence $\overline x^n=0$. Thus $\overline x$ is a nonzero nilpotent element in
grad$(R,I)$. This is a contradiction. Therefore for every positive integer $c$,
the ideal $I^c$ coincides with its integral closure in $R$.

\centerline{}

(3.3.2) LEMMA. Assume that $R$ is a normal noetherian domain and $I$ is a
nonzero nonunit normal ideal in $R$. Then the height zero primes of grad$(R,I)$ 
correspond in a one-to-one manner with the Rees valuations of $I$ as follows.
Upon letting $L$ be the quotient field of $R$ and $E(I)$ be the Rees ring of $I$
relative to $R$ with variable $Z$, $P\mapsto E(I)_P\cap L$ gives a 
bijection of the set of all minimal primes of $IE(I)$ onto the set of all
Rees valuation rings of $I$.

Observe that the height zero primes of grad$(R,I)$ form a nonempty finite set.
Also observe that under the natural epimorphism
$$
E(I)\to\frac{E(I)}{IE(I)}\to\text{grad}(R,I)
$$
the set of minimal primes of $IE(I)$ bijectively map onto the height zero primes
of grad$(R,I)$. Finally observe that, by (8.1)(VI) of \cite{AH1}, $E(I)$ is a
normal noetherian domain and $\mathfrak W(R,I)^{\mathfrak N}=\mathfrak W(R,I)$,
and hence the Rees valuation rings $\mathfrak D^*(R,I,I)$ coincide with the
dicritical divisors $\mathfrak D(R,I,I)$.

\centerline{}

PROOF. Follows from (6.11) of \cite{AH2}.

\centerline{}

(3.4) NOTE. Let $J=(y_1,\dots,y_b)R\subset(x_1,\dots,x_a)R=I$ be
any finitely generated ideals in a domain $R$ with quotient field $L$.
Then we have the following.

(I) $J$ is a reduction of $I\Leftrightarrow JV=IV$ for all 
$V\in\overline D(L/R)$.

(II) If $J$ is a reduction of $I$ then, for every $p\in\mathbb N$, the ideal 
$(y_1^p,\dots,y_b^p)R$ is a reduction of $I^p$. 

\centerline{}

PROOF OF (I). First suppose that $J$ is a reduction of $I$. Then
$JI^n=I^{n+1}$ for some $n\in\mathbb N$. Let there be given any 
$V\in\overline D(L/R)$. Since $J$ and $I$ are finitely generated, they extend to 
principal ideals $JV=\beta V$ and $IV=\alpha V$ with $\beta\in J$ and $\alpha\in I$.
Therefore, since $JI^n=I^{n+1}$, we get $\alpha\beta^n V=\beta^{n+1}V$. 
Consequently $\alpha V=\beta V$. Hence $JV=IV$.

Next suppose that $JV=IV$ for all $V\in\overline D(L/R)$.
Then $I$ is integral over $J$ by (3.6) below.
Therefore, by (4.4) of \cite{AH1}, for every $i\in\{1,\dots,a\}$ we can find
$n_i\in\mathbb N_+$ such that $x_i^{n_i}\in JI^{n_i-1}$.
Let $n=n_1+\dots+n_a$. We claim that then $I^{n+1}\subset JI^n$. From this it will
follow that $JI^n=I^{n+1}$ and hence $J$ is a reduction of $I$. To prove the claim
it suffices to show that $x_1^{m_1}\dots x_a^{m_a}\in JI^n$ whenever
$m_1+\dots+m_a=n+1$. But the last equation implies $m_i\ge n_i$ for some $i$
and hence $x_i^{m_i}\in JI^{m_i-1}$. 
Therefore $x_1^{m_1}\dots x_a^{m_a}\in JI^n$.

\centerline{}

PROOF OF (II). Assume that $J$ is a reduction of $I$. Then for all
$V\in\overline D(L/R)$, by (I) we have $JV=IV$ and hence $J^pV=I^pV$ for
every $p\in\mathbb N$. But clearly $J^pV=(y_1^p,\dots,y_b^p)V$ and hence
$(y_1^p,\dots,y_b^p)V=I^pV$. Therefore again by (I) we see that
$(y_1^p,\dots,y_b^p)R$ is a reduction of $I^p$.

\centerline{}

(3.5) ZARISKI'S THEOREM (see Theorem 1 on page 350 of Volume II of \cite{Zar} 
and the beginning of Section 8 of \cite{AH1}).
Let $\overline R$ be the integral closure of a domain $R$ in its quotient field
$L$. Then the completion of any ideal $J$ in $R$ coincides with the integral
closure of $J$ in $\overline R$.

\centerline{}

(3.6) ZARISKI'S COROLLARY. In the situation of (3.5), let $J\subset I$ be any ideals
in $R$ such that $JV=IV$ for all $V\in\overline D(L/R)$. Then $I$ is integral
over $J$.

\centerline{}

PROOF. By the definition of completion, the completion of $I$ coincides with the
completion of $J$. Therefore by (3.5) the integral closure of $I$ in 
$\overline R$ coincides with the integral closure of $J$ in $\overline R$.
Therefore $I$ is integral over $J$.

\centerline{}

(3.7) REMARK. Comparing the concrete (3.2) to the abstract (3.3), we note that in
(3.3) we are saying that the dicriticals are in a bijective correspondence with the
height zero primes of grad$(R,I)$ under the assumption that $R$ is a normal
noetherian domain and grad$(R,I)$ is reduced, without explicitly describing how to
find these height zero primes or how many such height zero primes there are. In
(3.2) we are taking a concrete hypersurface case and giving an explicit description
of the dictitical divisors of the maximal ideal by using affine QDTs.

\centerline{}

{{\bf Section 4: Contact Numbers.}} 
Given any prime divisors $V,W$ of a two dimensional regular local domain $R$, 
i.e., elements $V,W$ of $D(R)^\Delta$,
in Proposition (3.5) of Abhyankar's paper \cite{Ab11}, the
{\bf contact number} $c(R,V,W)$ of $V$ with $W$ at $R$ was defined by putting
$$
c(R,V,W)=V(\zeta_R(W))
$$
where we recall that for any $a\in\text{QF}(V)$ we have put 
$$
V(a)=\ord_V(a)
$$
and for any ideal $I$ in a noetherian subring $S$ of $V$ we have put
$$
V(I)=\text{min}\{V(a):a\in I\}
$$ 
with the understanding that $V(I)=\infty$ if $I\subset\{0\}$.

Our goal in this Section is to prove the following commutativity Theorem (4.6)
about contact numbers. Our main tool will be Lemma (6.11) of our previous paper 
\cite{AH2}. First we shall introduce some terminology and prove a string of Lemmas.

Recall that $H_T:T\to H(T)=T/M(T)$ denotes the residue class epimorphism
of any quasilocal ring $T$, and for any subring $S$ of $T$ let us put
$$
\chi(S,T)= [K':K]
$$
where $K'$ is the algebraic closure of $K=\text{QF}(H_T(S))$ in $H(T)$.

 For any subring $S$ of a ring $T$ we put
$$
\overline S^T=\text{the integral closure of $S$ in $T$}
$$
and we note that $\overline S^T$ is a subring of $T$,
and for any ideal $J$ in $S$ we put
$$
J^{-T}=\text{the integral closure of $J$ in $T$}
$$
and we note that $J^{-T}$ is an ideal in $\overline S^T$. 

Recall that the {\bf length} of a module $M$ over a ring $D$ is denoted by 
$$
\ell_D M.
$$
Note that $M$ is a {\bf finite $D$-module} means $M$ is finitely generated
as a module.

By $\iota(\mathfrak a,\mathfrak a';R)$ we denote the 
{\bf intersection multiplicity} of any nonzero principal ideals
$\mathfrak a,\mathfrak a'$ in a two dimensional regular local domain $R$;
for respective generators $a,a'$ of $\mathfrak a,\mathfrak a'$, we may write 
$\iota(a,a';R)$ instead of $\iota(aR,a'R;R)$; recall that 
$$
\iota(a,a';R)=\ell_R(R/(a,a')R)
$$
and note that this is zero or a positive integer or infinity 
according as the ideal $(a,a')R$ is the unit ideal or an $M(R)$-primary ideal 
or is contained in a nonzero nonunit principal ideal. 

Let $t,t^*$ be independent indeterminates over a field $L$. 
Referring to the beginning of Section 2 of \cite{Ab13} for the definition of
the $t$-extension $R^t$ of any subring $R$ of $L$, we define the 
{\bf $(t,t^*)$-extension} $R^{t,t^*}$ of $R$ by putting $R^{t,t^*}=(R^t)^{t^*}$. 
The ring $R^{t,t^*}$ can be directly defined as the localzation of 
$R[t,t^*]$ at the multiplicative set of all polynomials $f(t,t^*)\in R[t,t^*]$
whose coefficient generate the unit ideal in $R$.
Note that if $R$ is a regular local domain of dimension $d$ with quotient field
$L$ then $R^t$ and $R^{t,t^*}$ are regular local domains of dimension $d$ with
quotient fields $L^t=L(t)$ and $L^{t,t^*}=L(t,t^*)$ respectively.

\centerline{}

NOTE (4.1). Let $J$ be an $M(R)$-primary ideal in a two dimensional regular
local domain, and let $I=J^{-R}$. Then by (4.5) and (4.6) of \cite{AH1} we 
know that $I$ is an  $M(R)$-primary ideal in $R$ and, by (3.5), 
$I$ coincides with the completion of $J$.
By (8.1)(VI) of \cite{AH1} and Zariski's Theorem stated as ZQT in Section 2 of 
\cite{AH1}, $I$ is a normal ideal in $R$, 
$\mathfrak D(R,J)=\mathfrak D(R,I,I)=$ a set consisting of a finite number
of distinct DVRs with $h\in\mathbb N_+$, and we have
$$
I=\prod_{1\le i\le h}\zeta_R(V_i)^{n(i)}
\quad\text{ with }\quad n(i)\in\mathbb N_+.
$$

\centerline{}

DEFINITION (4.2). In the above set-up of (4.1), let $V^*_1,\dots,V^*_{h^*}$
be a finite number of distinct DVRs with $h^*\in\mathbb N_+$ and let
$$
I^*=\prod_{1\le i\le h^*}\zeta_R(V^*_i)^{n^*(i)}
\quad\text{ with }\quad n^*(i)\in\mathbb N_+.
$$
Let us define the {\bf contact number} $c(R,I,I^*)$ of $I$ with $I^*$ at $R$
by putting
$$
c(R,I,I^*)=\sum_{1\le i\le h,1\le j\le h^*}n(i)n^*(j)\chi(R,V_i)c(R,V_i,V^*_j).
$$
Note that this is always a positive integer.
In (4.6) we shall prove a commutativity property of $c(R,I,I^*)$ and, as a
special case, it will imply a commutativity property of $c(R,V,W)$.

\centerline{}

LEMMA (4.3). Let $D$ be a one dimensional local domain with quotient field $L$
and let $B=\overline D^L$. Assume that $B$ is finite $D$-module. Then $B$
is a Dedekind Domain having only a finite number of distinct nonzero prime 
ideals $P_1,\dots,P_h$ with $h\in\mathbb N_+$, and upon letting $V_i=B_{P_i}$ 
we have that $V_1,\dots,V_h$ are DVRs with $B=V_1\cap\dots\cap V_h$ and
$M(V_i)\cap D=M(D)$ for $1\le i\le h$.
Moreover, for any $0\ne z\in M(D)$ we have
$$
\ell_D(D/zD)=\sum_{1\le i\le h}e_if_i\;\text{ where }\; e_i=V_i(z) 
\;\text{ and }\; f_i=\chi(D,V_i).
$$

\centerline{}

PROOF. Everything except the ``Moreover'' part is well-known. Now for any 
$0\ne z\in D$ we 
$$
zB=\prod_{1\le i\le h} P_i^{e_i}\;\text{ with }\;e_i=V_i(z).
$$
Considering the natural surjective maps 
$$
\phi:B\to B/zB 
\;\text{ and }\;\phi_i:B\to B/P_i^{e_i}
$$
and
$$
\psi_i:B/P_1^{e_1}\oplus\dots\oplus B/P_h^{e_h}\to B/P_i^{e_i}
$$
of $B$-modules, by chinese remaindering we get a unique isomorphism
$$
\psi:B/zB\to B/P_1^{e_1}\oplus\dots\oplus B/P_h^{e_h}
$$
of $B$-modules such that for all $x\in B$ we have
$\psi_i(\psi(\phi(x))=\phi_i(x)$ for $1\le i\le h$. Clearly $\psi$ is also an
isomorphism of $D$-modiles, and hence
$$
\ell_D(D/zD)=\sum_{1\le i\le h}\ell_D(B/P_i^{e_i}).
$$
We shall show that for each $i$ we have $\ell_D(B/P_i^{e_i})=e_if_i$ and this
will complete the proof. Again by chinese remaindering we see that each $P_i$
is generated by a single nonzero element $y_i$.
We have a sequence
$$
P_i^{e_i}\subset P_i^{e_i-1}\subset P_i=P_i^0=B
$$ 
of $B$-modules and multiplication by $y_i^j$ gives a 
$B$-module isomorphism $B/P_i\to P_i^j/P_i^{j+1}$ and hence
a $D$-module isomorphism $B/P_i\to P_i^j/P_i^{j+1}$.  Therefore for 
$0\le j\le e_i-1$ we have $\ell_D(B/P_i)=\ell_D(P_i^j/P_i^{j+1})$. 
Consequently 
$$
\ell_D(D/P_i^{e_i})=e_i\ell_D(B/P_i).
$$
Now $V_i/M(V_i)$ and $B/P_i$ are isomorphic as $B$-modules, and hence also
as $D$-modules.  Therefore it suffices to show that
$\ell_D H(V_i)=f_i$. Clearly $H(V_i)$ is annihilated by
$M(D)$ and hence $H(V_i)$ may be regarded as a $D/M(D)$-module and we have 
$\ell_D H(V_i)=\ell_{D/M(D)}H(V_i)$. Obviously $\ell_{D/M(D)}H(V_i)=f_i$.

\centerline{}

LEMMA (4.4). Let $R$ be a normal noetherian domain with quotient field $L$. 
Let $I$ a nonzero nonunit normal ideal in $R$. Let $J$ be an ideal in $R$ such
that $I=J^{-R}$. Let $0\ne x\in I$ be such that for every
$V\in\mathfrak D(R,I,I)$ we have $V(x)=V(I)$. Let $C=R[Jx^{-1}]$ and
$A=R[Ix^{-1}]$. Then $A$ is a finite $C$-module with $A=\overline C^L$.

\centerline{}

PROOF. By (6.11) of \cite{AH2} we know that 
$\mathfrak W(R,I)^{\mathfrak N}=\mathfrak W(R,I)$ and, because $0\ne x\in I$,
clearly $\mathfrak V(A)$ is an affine piece of $\mathfrak W(R,I)$, i.e.,
$$
A=\bigcap_{S\in\mathfrak V(A)}S\quad\text{ with }\quad 
\mathfrak V(A)\subset\mathfrak W(R,I)
$$
and hence $\overline A^L=A$. Since $I$ is integral over $J$, any $y\in I$
satisfies an equation of the form 
$$
y^n+y_1 y^{n-1}+\dots+y_iy^{n-i}y_n=0\;\text{ where }\;
y_i\in J^i\;\text{ for }\;1\le i\le n
$$ 
with $n\in\mathbb N_+$. We can write $y_i$ as a finite sum 
$\sum y_{i1}\dots y_{ii}$ with $y_{i1},\dots,y_{ii}$ in $J$. Dividing the 
displyed equation by $x^n$ we get 
$$
\left(\frac{y}{x}\right)^n+\sum_{1\le i\le n}\left(\sum
\left(\frac{y_{i1}}{x}\right)\dots\left(\frac{y_{ii}}{x}\right)\right)
\left(\frac{y}{x}\right)^{n-i}=0
$$
with $\frac{y_{ij}}{x}\in Jx^{-1}$ for all $i,j$. 
Since $I$ is a finitely generated ideal, it follows that $A$ is a finitely
geneartaed ring extension of $C$ and every element of $A$ is integral over $C$.
Therefore by (E2) on page 161 of \cite{Ab8} we conclude that 
$A$ is a finite $C$-module with $A=\overline C^L$.

\centerline{}

NOTE (4.5). Let $R$ be a normal quasilocal domain with quotient field $L$, and
let $C=R[F/G]$ where $F,G$ are nonzero elements in $R$ such that $F/G\not\in R$
and $G/F\not\in R$. Let $Q=M(R)C$ and let $Z$ be an indeterminate.

(1) By the bracketed remark on pages 75-76 of \cite{Ab2}, there exists a unique 
ring epimorphism $\phi:C\to H(R)[Z]$ with $\phi(F/G)=Z$ such that for all
$a\in R$ we have $\phi(a)=H_R(a)$. By the said bracketed remark we also see
that ker$(\phi)=Q$ and hence $Q$ is a nonzero depth-one prime ideal in $C$.

(2) Let us observe that if $R$ is a two dimensional normal local domain then 
$C$ is a two dimensional noetherian domain and $Q$ is a height-one prime ideal 
in $C$. To see this note that, by Lemma (T30) on page 235 of \cite{Ab8}, $R[t]$ 
is a three dimensional noetherian domain, and by sending $t$ to $-F/G$ we get an 
$R$-epimorphism (= a ring epimorphism which keeps $R$ elementwise fixed) 
$\mu:R[t]\to C$ with $\mu(\Phi)=0$ where
$$
\Phi=F+tG
$$
and hence $C$ is a two dimensional noetherian domain and $Q$ is a height-one 
prime ideal in $C$. 

(3) Henceforth assume that $R$ is a two dimensional regular local domain and 
dividing $F$ and $G$ by their GCD (which is determined up to a unit in $R$)
let us arrange that the ideal $J=(F,G)R$ is $M(R)$-primary. Also let $I=J^{-R}$ 
and $A=R[I/G]$. By (2) we see that $Q$ is a height-one prime ideal in 
the two dimensional noetherian domain $C$, and hence upon letting $D=C_Q$ we see 
that $D$ is a one dimensional local domain and $\mu$ uniquely extends to an
$R$-epimorphism $\nu:R^t\to D$. Note that now 
$$
\text{ker}(\mu)=\Phi R[t]\quad\text{ and }\quad\text{ker}(\nu)=\Phi R^t.
$$
By (4.1) we know that $I$ is a normal ideal in $R$ 
and hence by (4.4) we see that $A$ is a finite $C$-module with $A=\overline C^L$. 
Consequently, upon letting $B=A_{C\setminus Q}$ we conclude that $B$ is a finite 
$D$-module with $B=\overline D^L$ and $B$ is a Dedekind Domain having only a 
finite number of distinct nonzero prime ideals $P_1,\dots,P_h$ with 
$h\in\mathbb N_+$, and upon letting $V_i=B_{P_i}$ 
for $1\le i\le h$ we have that $V_1,\dots,V_h$ are DVRs with 
$B=V_1\cap\dots\cap V_h$. In view of Lemma 8.3 of \cite{AH1}, 
by (4.1) it follows that 
$V_1,\dots,V_h$ are all the distinct members of $\mathfrak D(R,J)$ and we have
$$
I=\prod_{1\le i\le h}\zeta_R(V_i)^{n(i)}
\quad\text{ with }\quad n(i)\in\mathbb N_+.
$$
By (4.3) we see that for any $0\ne z\in M(D)$ we have
$$
\ell_D(D/zD)=\sum_{1\le i\le h}\chi(D,V_i)\,V_i(z).
$$
By the description given in Lemma 8.3(II) of \cite{AH1} we see that
for $1\le i\le h$ we have
$$
\chi(D,V_i)=n(i)\chi(R,V_i).
$$
Combining the above two displays we obtain
$$
\ell_D(D/zD)=\sum_{1\le i\le h}n(i)\chi(R,V_i)\,V_i(z).
$$

(4) In the set-up of (3) let 
$$
I^*=\prod_{1\le i\le h^*}\zeta_R(V^*_i)^{n^*(i)}
\quad\text{ with }\quad n^*(i)\in\mathbb N_+
$$
where $V^*_1,\dots,V^*_{h^*}$ are pairwise distinct members of $D(R)^\Delta$ 
with $h^*\in\mathbb N_+$. Also let there be given any
$\Phi^*\in I^*R[t]$ such that $(\Phi,\Phi^*)R^t$ is $M(R^t)$-primary and 
$$
(V_i)^t(\Phi^*)=V_i(I^*)=V_i(\mu(\Phi^*))\quad\text{ for }\quad 1\le i\le h.
\leqno(\sharp)
$$
We CLAIM that then
$$
\iota(\Phi,\Phi^*;R^t)=c(R,I,I^*).
$$

\centerline{}

PROOF OF THE CLAIM. Upon letting $z=\mu(\Phi^*)$ we get 
a nonzero element $z$ in $M(D)$ such that $R^t/(\Phi,\Phi^*)R^t$ and  
$D/zD$ are isomorphic $R^t$-modules and we have
$$
\ell_{R^t}(R^t/(\Phi,\Phi^*)R^t)=\ell_D(D/zD).
$$
By definition the LHS equals $\iota(\Phi,\Phi^*;R^t)$ and by (3) the RHS
equals
$$
\sum_{1\le i\le h}n(i)\chi(R,V_i)\,V_i(z)
$$
and hence we get
$$
\iota(\Phi,\Phi^*;R^t)=\sum_{1\le i\le h}n(i)\chi(R,V_i)\,V_i(z).
$$
Therefore by $(\sharp)$ we obtain
$$
\iota(\Phi,\Phi^*;R^t)=\sum_{1\le i\le h}n(i)\chi(R,V_i)\,V_i(I^*).
$$
But for $1\le i\le h$ we have
$$
V_i(I^*)=\sum_{1\le j\le h^*}n^*(j)V_i(\zeta_R(V^*_j)).
$$
By the above two equations we get
$$
\iota(\Phi,\Phi^*;R^t)=\sum_{1\le i\le h,1\le j\le h^*}
n(i)n^*(j)\chi(R,V_i)\,V_i(\zeta_R(V^*_j)).
$$
Therefore by the definition of $c(R,I,I^*)$ we conclude that
$$
\iota(\Phi,\Phi^*;R^t)=c(R,I,I^*).
$$

\centerline{}

THEOREM ON COMMUTATIVITY OF CONTACT NUMBERS (4.6). Let $R$ be a two dimensional
regular local domain with quotient field $L$. Let
$$
I=\prod_{1\le i\le h}\zeta_R(V_i)^{n(i)}
$$
where $h,n(1),\dots,n(h)$ are positive integers, and $V_1,\dots,V_h$ are pairwise
distinct members of $D(R)^\Delta$. Let
$$
I^*=\prod_{1\le i\le h^*}\zeta_R(V^*_i)^{n^*(i)}
$$
where $h^*,n^*(1),\dots,n^*(h^*)$ are positive integers, and 
$V^*_1,\dots,V^*_{h^*}$ are pairwise distinct members of $D(R)^\Delta$. 
Let $F,G,F^*,G^*$ be nonzero members of $M(R)$ such that
$$
((F,G)R)^{-R}=I\;\text{ and }\; ((F^*,G^*)R)^{-R}=I^*.
$$
Then we have the following.

(4.6.1) Upon letting
$$
\Phi=F+tG\;\text{ and }\;\Phi^*=F^*+t G^*
$$
and, assuming $(\bullet)$ $(\Phi,\Phi^*)R^t$ to be $M(R^t)$-primary, we have
$$
c(R,I,I^*)=\iota(\Phi,\Phi^*;R^{t})=c(R,I^*,I).
$$

(4.6.2) We always have
$$
c(R,I,I^*)=c(R,I^*,I).
$$

(4.6.3) Upon letting
$$
\Phi=F+tG\;\text{ and }\;\Phi^*=F^*+t^* G^*
$$
we have
$$
c(R,I,I^*)=\iota(\Phi,\Phi^*;R^{t,t^*})=c(R,I^*,I).
$$

\centerline{}

PROOF OF (4.6.1). By symmetry, it is enough to prove that
$$
\iota(\Phi,\Phi^*;R^{t})=c(R,I,I^*).
$$
Therefore, in view of the CLAIM of (4.5)(4), we only need to show that 
for $1\le i\le h$ we have
$$
(V_i)^t(\Phi^*)=V_i(I^*)=V_i(\mu(\Phi^*)).
$$
Clearly
$$
(V_i)^t(\Phi^*)=\text{min}(V_i(F^*),V_i(G^*))=V_i(I^*)
$$
and hence $(V_i)^t(\Phi^*)=V_i(I^*)$. Now
$$
\mu(\Phi^*)=F^*-(F/G)G^*
$$
and by Section 2 of \cite{Ab11} we know that $F/G$ is residually transcendental
over $R$ at $V_i$. Therefore $V_i(\mu(\Phi^*))=\text{min}(V_i(F^*),V_i(G^*))$.

\centerline{}

PROOF OF (4.6.2). Let $\Phi=F+tG$. Let
$$
\Phi^*=\begin{cases}
F^*+tG^*&\text{if $(\Phi,F^*+tG^*)R^t$ is $M(R^t)-primary$}\\
G^*+tF^*&\text{otherwise.}\\
\end{cases}
$$
Then $(\Phi,\Phi^*)R^t$ is $M(R^t)$-primary, and hence we are done by (4.5.1).

\centerline{}

PROOF OF (4.6.3). (See (4.3) on page 334 of \cite{Hun}).
By symmetry, it is enough to prove that
$$
\iota(\Phi,\Phi^*;R^{t,t^*})=c(R,I,I^*).
\leqno({}^{*})
$$
Let $(\overline R,\overline L)=(R^{t^*},L^{t^*})$ and
$(\overline I,\overline I^*)=(I\overline R,I^*\overline R)$ and
$$
(\overline V_1,\dots,\overline V_h,\overline V^*_1,\dots,\overline V^*_{h^*})
=((V_1)^{t^*},\dots,(V_h)^{t^*},(V^*_1)^{t^*},\dots,(V^*_{h^*})^{t^*}).
$$
Then $\overline V_1,\dots,\overline V_h$ are pairwise distinct members of 
$D(\overline R)^\Delta$ and $F,G$ are nonzero members of $M(\overline R)$ 
such that
$$
((F,G)\overline R)^{-\overline R}=\overline I
=\prod_{1\le i\le h}\zeta_{\overline R}(\overline V_i)^{n(i)}.
$$
Likewise, $\overline V^*_1,\dots,\overline V^*_h$ are pairwise distinct 
members of $D(\overline R)^\Delta$ and $F^*,G^*$ are nonzero members of 
$M(\overline R)$ such that
$$
((F^*,G^*)\overline R)^{-\overline R}=\overline I^*
=\prod_{1\le i\le h^*}\zeta_{\overline R}(\overline V^*_i)^{n^*(i)}.
$$
Clearly the ideal $(\Phi,\Phi^*)\overline R^t$ is 
$M(\overline R^t)$-primary and for $1\le i\le h$ we have 
$$
(\overline V_i)^t(\Phi^*)=\text{min}(\overline V_i(F^*),\overline V_i(G^*))
=\overline V_i(I^*)
$$
and hence $(\overline V_i)^t(\Phi^*)=\overline V_i(I^*)$. Also
$$
\mu(\Phi^*)=F^*+t^*G^*
$$
where $\mu:\overline R[t]\to C=\overline R[F/G]$ is the unique 
$\overline R$-epimorphism such that $\mu(t)=-F/G$ and 
ker$(\mu)=\Phi\overline R[t]$; therefore 
$\overline V_i(\mu(\Phi^*))=\text{min}(\overline V_i(F^*),\overline V_i(G^*))$.
It follows that
$$
(\overline V_i)^t(\Phi^*)=\overline V_i(\overline I^*)
=\overline V_i(\mu(\Phi^*))\quad\text{ for }\quad 1\le i\le h.
\leqno(\sharp)
$$
Hence, by putting a bar on the relevant quantities in the CLAIM of (4.5)(4), 
we conclude that
$$
\iota(\Phi,\Phi^*;\overline R^t)=c(\overline R,\overline I,\overline I^*).
$$
Clearly the above LHS equals the LHS of (*), and the above RHS equals the
RHS of (*). This establishes (*) and concludes the proof.

\centerline{}

NOTE (4.7). As noted in the proof of (4.6.2), the condition $(\bullet)$ of 
(4.6.1) is satisfied  by at least one of the two pairs
$(\Phi,F^*+tG^*)$ or $(\Phi,G^*+tF^*)$. It is also satisfied by letting
$\Phi=F+tG$ and $\Phi^*=F+\tau F$ where $\tau=$ any unit in $R$ excluding 
at most one value.

\centerline{}

NOTE (4.8).  By a {\bf testing curve} for $V\in D(R)^\Delta$, where $R$ is
a two dimensional regular local domain, we mean an element $\delta_V\in R$ which 
does ``something'' useful for $V$ at $R$. For instance 
we ask whether, for every $V\in D(R)^\Delta$, there is a testing curve 
$\delta_V$ such that $W(\delta_V)=W(\zeta_R(V))$ for all $W\in D(R)^\Delta$. 
To answer this negatively, let 
$$
V=o(R)
$$
and note that
$$
\zeta_R(V)=M\;\;\text{ with }\;\; V(\zeta_R(V))=1
$$
and for any $\theta\in R$ we have
$$
V(\theta)=V(\zeta_R(V))\Leftrightarrow\theta\in M\setminus M^2.
\leqno(1)
$$
Now let $(x,y)$ be generators of $M$, let $K$ be a coefficient set of $R$ and, 
for every $t\in K\cup\{\infty\}$, let
$$
(p_t,q_t)=\begin{cases}(y+tx,x)&\text{if }t\in K\\
(x,y)&\text{if }t=\infty\\ \end{cases}
$$
and
$$
\Gamma_t=\{\theta\in M\setminus M^2:(\theta R)+M^2=(p_t R)+M^2\}
$$
and
$$
S_t=R[p_t/q_t]_{(p_t/q_t,q_t)R[p_t/q_t]}\in Q_1(R)
$$
and 
$$
W_t=o(S_t)\in D(R)^\Delta.
$$ 
Note that then
$$
\begin{cases}
\text{$t\mapsto S_t$ gives a bijection of $K\cup\{\infty\}$}\\
\text{onto the set of all those members of $Q_1(R)$}\\
\text{which are residually rational over $R$}
\end{cases}
$$
and
$$
M\setminus M^2=\coprod_{t\in K\cup\{\infty\}}\Gamma_t=
\text{a partition into disjoint nonempty subsets.}
\leqno(2)
$$ 
Also note that for each $t\in K\cup\{\infty\}$ we have
$$
(p_t,q_t)R=M
$$
with
$$
W_t(\theta)=W_t(\zeta_R(V))=1\text{ for all }
\theta\in(M/\setminus M^2)\setminus\Gamma_t
\leqno(3)
$$ 
and
$$
W_t(\theta)=2\text{ for all }
\theta\in\Gamma_t.
\leqno(4)
$$ 
By (1) to (4) we see that
$$
\begin{cases}
\text{there is no $\delta_V\in R$}\\
\text{such that $V(\delta_V)=V(\zeta_R(V))$}\\
\text{and $W_t(\delta_V)=W_t(\zeta_R(V))$ for all $t\in\ K\cup\{\infty\}$}\\
\end{cases}
\leqno(5)
$$ 
but
$$
\begin{cases}
\text{for every $\tau\in K\cup\{\infty\}$ and every 
$\delta_V\in\Gamma_\tau$}\\
\text{we have $V(\delta_V)=V(\zeta_R(V))$}\\
\text{and $W_t(\delta_V)=W_t(\zeta_R(V))$ for all 
$t\in(K\cup\{\infty\})\setminus\{\tau\}$.}\\
\end{cases}
\leqno(6)
$$ 

\centerline{}

In \cite{AbA} it will be shown how the difficulty in definiing a good testing 
curve is removed by introducing the concept of ``free points'' of $V$.

\centerline{}

NOTE (4.9). Here is a {\bf Question} related to the definitions of
dicritical divisors and Rees valuation rings introduced in Section 2. 
Let $J$ be a nonzero nonunit ideal in a normal noetherian domain $R$ with quotient
field $L$, and for $\ne G\in J$ let $C=R[J/G]$ and $A=\overline C^L$. Then does 
every height-one prime of $GA$ contract to a height-one prime of $GC$?

\centerline{}

NOTE (4.10). Here is another {\bf question}. 
In (4.5)(4), we found a sufficient condition for the intersection
formula to work. That condition assumes an indeterminate inside $\Phi$
but not inside $\Phi^*$. The question is whether we can find a condition in 
which both $\Phi$ and $\Phi^*$ stay inside $R$. In other words, can we find a
sufficient (and necessary?) condition on $\Phi\in I$ and $\Phi^*\in I^*$ which 
would imply $\iota(\Phi,\Phi^*;R)=c(R,I,I^*)$?

\centerline{}

NOTE (4.11). Here is yet another {\bf question}. 
In (46) of Section 3, can you prove the normality of $R$ under more general
conditions?

\end{document}